\documentclass{article}
\usepackage[utf8]{inputenc} 
\usepackage[T1]{fontenc}    
\usepackage{lmodern}
\usepackage[margin=1in]{geometry}

\usepackage{natbib}
\setcitestyle{authoryear,round,citesep={;},aysep={,},yysep={;}}
\renewcommand{\cite}[1]{\citep{#1}}

\usepackage{amsthm} 
\usepackage[vvarbb]{newtxmath}
\usepackage{microtype}      

\usepackage{nicefrac}
\usepackage{graphicx,float,wrapfig}

\usepackage{mathtools}

\usepackage{enumitem}
\setitemize{noitemsep,topsep=0pt,parsep=0pt,partopsep=0pt}
\newlist{enumerate*}{enumerate*}{1}
\setlist[enumerate*]{label=(\arabic*), itemjoin={{, }},
  itemjoin*={{, and }}, after={.}}

\usepackage[ruled, linesnumbered]{algorithm2e} 

\usepackage[colorlinks,urlcolor=blue,citecolor=blue,linkcolor=blue]{hyperref}       
\usepackage{authblk}

\theoremstyle{plain} \numberwithin{equation}{section}
\newtheorem{theorem}{Theorem}[section]
\numberwithin{theorem}{section}

\theoremstyle{definition}

\newtheorem{observation}[theorem]{Observation}

\theoremstyle{plain}

\DeclareMathOperator*{\argmin}{argmin}

\newcommand{\setb}[1]{\left\{ #1 \right\}}
\newcommand{\R}{\mathbb R}
\newcommand{\norm}[1]{\left\| #1 \right\|}

\hypersetup{
  pdftitle={Dual Prices for Frank–Wolfe Algorithms (note)},
  pdfauthor={Gábor Braun, Sebastian Pokutta},
  pdfkeywords={Frank–Wolfe algorithm, conditional gradient,
    gradient descent, linear optimization oracle, convex optimization},
  pdfsubject={MSC2000 Primary 68Q32; 
    Secondary
    90C52 
  }}

\author[1,2]{Gábor Braun}
\author[1,2]{Sebastian Pokutta}
\affil[1]{Technische Universität Berlin, Germany}
\affil[2]{Zuse Institute Berlin, Germany\\ \texttt{\{braun,pokutta\}@zib.de}}

\title{Dual Prices for Frank–Wolfe Algorithms%
  \thanks{Research in this paper was partially supported by
    Research Campus MODAL funded by German Federal Ministry of
    Education and Research under grant 05M14ZAM.}
  \\ —a note—}
\date{January 6, 2021}

\begin{document}

\maketitle

\begin{abstract}
  In this note we observe that for constrained
  convex minimization problems \(\min_{x \in P}f(x)\)
  over a polytope \(P\),
  dual prices for the linear program \(\min_{z \in P} \nabla f(x) z\)
  obtained from linearization
  at approximately optimal solutions \(x\) have a similar interpretation
  of rate of change in optimal value
  as for linear programming,
  providing a convex form of
  sensitivity analysis.
  This is of particular interest for
  Frank–Wolfe algorithms (also called conditional gradients),
  forming an important class of first-order methods,
  where a basic building block is linear minimization of gradients of
  \(f\) over \(P\),
  which in most implementations already
  compute the dual prices as a by-product.
\end{abstract}

\section{Introduction}

We consider the constrained convex minimization problem
\begin{equation}
  \label{eq:min}
  \tag{minProb}
\min_{x \in P} f(x),
\end{equation}
where \(f\) is a smooth convex function and \(P\) is a compact convex
feasible region.
Our primary interest is \emph{first-order} algorithms
where access to \(f\) is provided by computing gradients
\(\nabla f(x)\) and function values \(f(x)\)
at any feasible point \(x\).  An important class
of first-order methods is formed by \emph{conditional gradient
  algorithms (also known as Frank–Wolfe algorithms)}, which access the
feasible region \(P\) solely through a \emph{linear minimization
  oracle}, i.e., presented with a linear objective \(c\) the
oracle returns \(\argmin_{x \in P} c \cdot x\).  This class
of algorithms has several advantages two of which are
\begin{enumerate*}
\item \emph{Projection-freeness:}  no projection to the domain \(P\)
  is needed
\item \emph{Sparsity:} iterates are represented as
  convex combination of a small number of vertices,
  usually at most one vertex per iteration
\end{enumerate*}
As an example, the simplest algorithm, namely,
the (vanilla) Frank–Wolfe Algorithm \cite{fw56,polyak66cg}
is recalled in Algorithm~\ref{alg:fw}, which however will not be used
in the rest of the paper.

\begin{wrapfigure}{I}{.5\textwidth}
\begin{algorithm}[H]
\caption{\label{alg:fw} Frank–Wolfe Algorithm (FW)}
$x_0\in P$ arbitrary\;
\For{$t=0$ \KwTo \dots}{
  $v_t \leftarrow \argmin_{z \in P} \nabla f(x_t) z$\\
  \(\gamma_{t} \leftarrow \argmin_{0 \leq \gamma\leq 1}
  f(x_{t} + \gamma (v_{t} - x_{t}))\)\\
  $x_{t+1} \leftarrow x_t + \gamma_t (v_t - x_t)$}
\end{algorithm}
\end{wrapfigure}

Many implementations of a linear minimization oracle already compute
dual prices for an optimal solution,
to verify optimality.
Therefore it is of interest to make use of this extra information.
The role of conditional gradient algorithms in this paper is only
as a practical example: making available dual prices at no
extra cost in most implementations.

We will provide an interpretation for dual prices similar to
sensitivity analysis for linear optimization:
dual prices at an optimal solution
\(x\) are the rate of change in optimal value under small changes to
the right-hand side \(b\) of constraints over a domain
\(\setb{z : A z \leq b}\) defined by linear inequalities.
We shall see this interpretation holds even for approximately optimal
solutions \(x\) while retaining
the additive error of the accuracy of \(x\).
Thus dual prices can then be used as customary,
e.g., in sensitivity analysis, to compute risk-free state
probabilities, (economic) shadow prices in e.g., energy systems, etc.

\subsection{Preliminaries}
\label{sec:preliminaries}

We briefly recall basic definitions.
Recall that a differentiable function \(f \colon P \to \mathbb{R}\)
is \emph{convex} if
for all $x,y\in P$
\begin{equation}
  \label{eq:convex}
  f(y) - f(x) \geq \nabla f(x)(y-x),
\end{equation}
We denote by \(\nabla f(x)\) the gradient of \(f\) as a row vector,
while \(P\) is considered to be contained in
a vector space of column vectors.
This will be convenient for dual prices, but we note that
the formalism is inherently free of coordinate choices
in the space of \(P\)
(treating \(\nabla f(x)\) as a linear function).

Further the
function \(f\) is \emph{$L$-smooth} if for all \(x,y \in P\)
\begin{equation}
  \label{eq:smooth}
  f(y) - f(x) \leq \nabla f(x) (y-x) + \frac{L}{2} \norm{y-x}^{2}.
\end{equation}
Here \(\norm{.}\) is any norm on the vector space containing \(P\) and
the value of \(L\) depends on the norm \(\norm{.}\).

Conditional gradient algorithms often use the \emph{Frank–Wolfe gap}
defined as:
\begin{equation}
  \max_{z \in P} \nabla f(x) (x-z) = \nabla f(x) (x - v),
\end{equation}
with \(v\) a point of \(P\) minimizing
\(\min_{z \in P} \nabla f(x) z\).
In the context of Frank–Wolfe algorithms, these minimizers are
called \emph{Frank–Wolfe vertices} at \(x\)
(even though not all are vertices, but only vertex minimizers are
used in practice).   They are usually used to define the
next iterate as in the vanilla variant in
Algorithm~\ref{alg:fw}.
Recall that by convexity
\begin{equation}
  0 \leq f(x) - f(x^*) \leq \nabla f(x)(x - x^{*})
  \leq \nabla f(x) (x - v)
  .
\end{equation}
Here and below \(x^{*}\) is an optimal solution to
\(\min_{x \in P} f(x)\).
Thus the Frank–Wolfe gap is an upper bound to
the primal gap \(f(x) - f(x^{*})\) and it is \(0\)
at optimal solutions to \(\min_{z \in P} f(z)\).
As such the Frank–Wolfe gap is useful as a proxy for the primal gap,
while it also provides a lower bound to the optimal value:
\(f(x) - \nabla f(x) (x - v) \leq  f(x) - \nabla f(x) (x - x^{*})
\leq f(x^*)\).

\section{Dual Prices in convex minimization}
\label{sec:dual-prices-context}

We recall dual prices for linear optimization applied to our context here.  Let
\(P = \setb{z : A z \leq b}\) with
\(A \in \R^{m \times n}, b \in \R^n\) be a polytope and let \(f\) be a
convex function and let \(x, v \in P\) be arbitrary. Recall that
\emph{strong duality} states that \(v \in P\) is a minimizer
for the linear program
\(\min_{z: Az \leq b} \nabla f(x) z\), i.e.,
\(v = \argmin_{z: Az \leq b} \nabla f(x) z\) if and only if
there is a nonnegative combination of constraints certifying optimality,
i.e., a vector \(0 \leq \lambda \in \R^m\) whose entries are
multipliers called \emph{dual prices} satisfying
\begin{subequations}
  \label{eq:grad-strong-duality}
  \begin{align}
    \label{eq:gradDecomp}
    \nabla f(x) &= - \lambda A, & \lambda &\geq 0,
    \\
    \label{eq:fvalDecomp} \nabla f(x) v = \min_{z: Az \leq b}
    \nabla f(x) z &= - \lambda b.
  \end{align}
\end{subequations}
The second equality can be replaced with \emph{complementary
  slackness}, namely, that \(v\) satisfies with equality the
constraints of \(A z \leq b\) whose multiplier in \(\lambda\) is
positive (i.e., \(a_{i} v = b_{i}\) for \(\lambda_{i} > 0\),
where \(a_{i}\) is row \(i\) of \(A\), and \(\lambda_{i}\) and
\(b_{i}\) are entry \(i\) of \(\lambda\) and \(b\) respectively),
or short
\(\lambda (b - A v) = 0\).
A \emph{primal-dual pair} is a pair \((v, \lambda)\)
satisfying the strong duality conditions stated in
Equation system~\eqref{eq:grad-strong-duality}.
Obviously \(v\) is a Frank–Wolfe vertex \(v\) at $\nabla f(x)$ and the
\emph{Frank–Wolfe gap at \(x\) equals the complementarity gap}, i.e.,
\begin{equation}
  \label{eq:4}
  \nabla f(x) (x-v) = - \lambda A (x-v) = \lambda (b - A x).
\end{equation}

Common implementations of a linear optimization oracle
naturally compute dual prices for a Frank–Wolfe vertex \(v\),
which is indeed a vertex of \(P\).
For example, the widely used simplex algorithm
internally operates with data providing
a candidate \((v, \lambda\)) for a primal-dual pair,
where \(v\) is a vertex of \(P\) but \(\lambda\) may violate the
nonnegativity condition, which is then incrementally improved
to a primal-dual pair.

Recall that the celebrated \emph{Slater's condition of optimality} (a special
case of the \emph{Karush–Kuhn–Tucker condition} for convex functions) is the
strong duality form of the optimality condition for \(\min_{z: A z \leq b}
f(z)\):
a point \(x\) is an optimal solution to \(\min_{z: A z \leq b} f(z)\)
if and
only if \(x\) is an optimal solution to \(\min_{z: Az \leq b} \nabla f(x) z\),
i.e., \((x,\lambda)\) is a primal-dual pair for the linear program \(\min_{z: Az
  \leq b} \nabla f(x) z\) for some \(\lambda\); equivalently,
there are dual prices $\lambda$ for \(x\)
under the linear objective $\nabla f(x)$.

In practical implementations, e.g., of Frank–Wolfe algorithms, this
means that dual prices \(\lambda\) for the optimal solution \(x^*\) can
be obtained as dual prices for the Frank–Wolfe vertex \(v\) associated
with \(\nabla f(x^*)\).

\subsection{Dual Prices for approximately optimal solutions:
Sensitivity in \(b\)}

In practice we rarely have exact optimal solutions to convex
minimization problems (even within the limit of numerical accuracy),
and we are usually satisfied with a good approximate solution
with, e.g., an additive error in function
value of at most \(\varepsilon\).
For Frank–Wolfe algorithms the usual stopping criterion is
an upper bound on the Frank–Wolfe gap (sometimes also called dual gap)
\(\max_{z: Az \leq b} \nabla f(x) (x - z) \leq \varepsilon\) as
due to the linear minimizations the Frank–Wolfe gap is essentially
computed anyway.  As such we will now
consider the case of approximately optimal solutions.
To this end let
\(v = \argmin_{z: Az \leq b} \nabla f(x) z\) be the
Frank–Wolfe vertex at \(x\) and let \(0 \leq \lambda \in \R^m\) be
associated dual prices as
in Equation system~\eqref{eq:grad-strong-duality}
above.

A common interpretation of \(\lambda\)
in the context of linear programs is the
rate of change in optimal value as a function of change to the
constant term \(b\), i.e.,
\(\min_{z: Az \leq b'} \nabla f(x) z = \lambda b'\) for \(b'\)
close to \(b\); if \(\lambda\) is not a unique dual solution it needs
to be chosen depending on \(b'\).  Morally, the optimal value changes
by \(\lambda (b - b')\), while the dual solution \(\lambda\) does not
change.

The next observation carries this sensitivity analysis over to smooth
convex functions and approximately optimal solutions: in this case we
will incur additional error terms due to
\begin{enumerate*}
\item
  non-linearity of the objective function
\item approximate optimality (with no error for optimal solutions)
\end{enumerate*}
There are many common assumptions on
the objective convex function \(f\)
to bound its non-linearity.
For the sake of exposition, we assume the most common one,
namely, smoothness,
which is only needed for the last inequality
in Equations~\eqref{eq:2} and \eqref{eq:3}.
For other assumptions the error term
\(\nicefrac{L \norm{v' - v}^{2}}{2}\)
should be replaced accordingly.
\begin{observation}
  Let \(f\) be an \(L\)-smooth convex function
  over a convex domain containing
  the polytopes \(P = \{ z : A z \leq b\}\)
  and \(P' = \{ z : A z \leq b'\}\).
  Let \(x \in P\) and \(v = \argmin_{z \in P} \nabla f(x) z\).
  Similarly, let \(v' = \argmin_{z \in P'} \nabla f(x) z\).
  Assume that \(x' \coloneqq x - v + v' \in P'\).  Then
  \begin{align}
    \label{eq:1}
    f(x) - \nabla f(x) (x - v)
    &\leq
    \min_{z \in P} f(z)
    \leq
    f(x)
    \\
    \label{eq:2}
    f(x) - \nabla f(x) (x - v) + \nabla f(x) (v' - v)
    &\leq
    \min_{z \in P'} f(z)
    \leq
    f(x')
    \leq
    f(x) + \nabla f(x) (v' - v)
    + \frac{L}{2} \norm{v' - v}^{2}
    \intertext{When \(\lambda\) is a common dual solution for both
      \(v\) in \(P\) and \(v'\) in \(P'\)}
    \label{eq:3}
    f(x) - \nabla f(x) (x - v) + \lambda (b - b')
    &\leq
    \min_{z \in P'} f(z)
    \leq
    f(x')
    \leq
    f(x) + \lambda (b - b')
    + \frac{L}{2} \norm{v' - v}^{2}
    .
  \end{align}
\end{observation}
To justify the assumption \(x' \in P'\),
we note that it holds
when \(b'\) is sufficiently close to \(b\) and
\(x\) is sufficiently close to the optimal solution \(x^{*}\)
(with \(v\) and \(v'\) appropriately chosen depending on \(b'\)).
Intuitively, a neighborhood of \(v'\) in \(P'\) is just a translation
of a neighborhood of \(v\) in \(P\) and \(x\) is well inside the
neighborhood to be preserved by translation.
Let us split the defining linear inequalities \(A z \leq b\) for \(P\)
into two:
let \(A_{=} z \leq b_{=}\) be the subsystem which \(x\) satisfies with
equality (describing the boundary of the neighborhood at \(v\)),
and \(A_{<} z \leq b_{<}\) be the subsystem with
inequalities which \(x\) satisfies with strict inequality
(describing far away parts of \(P\)),
i.e., \(A_{=} x = b_{=}\) and \(A_{<} x < b_{<}\).
We claim that when \(x\) is close enough to \(x^{*}\) then
\(A_{=} v = b_{=}\) (regardless of the choice of \(v\)).
In geometrical terms the claim means that \(v\) is contained in
the minimal face containing \(x\).
To verify it,
let \(F\) be the face of \(P\) containing \(x^{*}\) in its relative
interior (allowing \(F = P\)), which is a minimal solution to
\(\min_{z \in P} \nabla f(x^{*}) z\)
by optimality.
When \(x\) is close to \(x^{*}\) then
\begin{enumerate*}
\item
  all minimal solutions to \(\min_{z \in P} \nabla f(x) z\)
  lie in \(F\), too
\item
  every hyperface of \(P\) containing \(x\)
  also contains \(x^{*}\) and hence \(F\)
  is contained in the minimal face containing \(x\)
\end{enumerate*}

As for linear programs, if \(b'\) is
sufficiently close to \(b\) then for some choice of optimal solutions
\(v\) and \(v'\) (recall they need not be unique)
they have a common dual solution
\(\lambda\) and \(v'\) is sufficiently close to \(v\).

After these preliminaries, we verify \(x' \in P'\).
First we deal with the inequalities \(x\) satisfy with equality for
\(P\), which turns out to be the easy case:
\(A_{=} x' = A_{=} x - A_{=} v + A_{=} v'
\leq b_{=} - b_{=} + b'_{=} = b'_{=}\).
For the other inequalities note that
\(b'_{<} - A_{<} x' = (b_{<} - A_{<} x)
+ (b'_{<} - b_{<}) - A_{<} (v' -  v)\).
As \(b_{<} - A_{<} x\) is strictly positive, with \(b'\) close enough
to \(b\) (and hence \(v'\) close enough to \(v\)),
the other terms on the right-hand side are small enough for the
right-hand side remaining positive, i.e., \(A_{<} x' < b'_{<}\).

\begin{proof}
The bounds for the polytope \(P\) in \eqref{eq:1} are well-known and presented for
comparison only; they easily follow from convexity and minimality of \(v\).

Equation~\eqref{eq:2} provides
the same bounds for the polytope \(P'\),
even though
the points at which we compute the lower and the upper bound
might differ.
The left-hand side of the first inequality is
\(f(x) - \nabla f(x) (x-v')\) (written in a form to ease comparison
with Equation~\eqref{eq:1}).
While \(x\) might not be contained in \(P'\), it does not affect the
validity of the inequality.
The second inequality explicitly uses the assumption \(x' \in P'\),
and the last inequality is just the smoothness inequality for \(f\),
using \(x' - x = v' - v\).

Finally, observe that \(\lambda (b - b') =  \nabla f(x) (v' - v)\) by
the definition of dual prices, leading to Equation~\eqref{eq:3}.
\end{proof}

\bibliographystyle{abbrvnat}
\bibliography{bibliography}

\end{document}